\documentclass[10pt,letterpaper]{article}

\makeatletter

\usepackage{latexsym}

\makeatother

\usepackage{blindtext,titlefoot}
\usepackage{authblk}
\usepackage{marvosym}
\usepackage[utf8]{inputenc}
\usepackage[english]{babel}
\usepackage{amsmath}
\usepackage{amsfonts}
\usepackage{amssymb}
\usepackage{makeidx}
\usepackage{graphicx}
\usepackage{lmodern}
\usepackage{kpfonts}
\usepackage{dsfont}
\usepackage{cancel}
\usepackage{amsthm} 
\usepackage[left=2cm,right=2cm,top=2cm,bottom=2cm]{geometry}
\usepackage{subcaption}
\usepackage{multirow}
\usepackage{float}
\usepackage{url}
\usepackage{breakurl}
\usepackage[breaklinks]{hyperref}
\usepackage{multicol}
\usepackage{mathtools, cuted}
\usepackage{lipsum}
\usepackage{booktabs}
\usepackage{comment}
\usepackage{cite}

\usepackage{fancyhdr}
\fancypagestyle{firstpage}{\fancyhf{}
\setcounter{page}{1}
\rfoot{\thepage}
}

\providecommand{\nset}[1]{
\mathbb{#1}
}
\providecommand{\set}[1]{
\left\{#1\right\}
}
\providecommand{\com}[1]{``#1"}

\providecommand{\ifr}[5]{
{}^{#1}_{#2}{#3}_{#4}^{#5}
}
\providecommand{\gam}[1]{
\Gamma\left(#1 \right)
}

\providecommand{\norm}[1]{
\left\lVert #1 \right\rVert
}
\providecommand{\abs}[1]{
\left\lvert #1 \right\rvert
}
\providecommand{\ds}[1]{
\displaystyle #1
}
\providecommand{\der}[3]{
\dfrac{#1^{#3} }{ #1 #2^{#3}}
}

\providecommand{\re}[1]{
\hbox{Re}\left(#1 \right)
}
\providecommand{\im}[1]{
\hbox{Im}\left(#1 \right)
}
\providecommand{\rnd}[2]{
\hbox{Rnd}_#2\left(#1\right)
}

\newtheorem{theorem}{ Theorem}[section]
\newtheorem{definition}[theorem]{Definition}
\newtheorem{proposition}[theorem]{Proposition}
\newtheorem{corollary}[theorem]{Corollary}
\newtheorem{example}[theorem]{Example}

\usepackage{enumitem} 
\setlist[itemize]{noitemsep} 

\usepackage{titlesec} 
\titleformat{\section}[block]{\large\bfseries\scshape\centering}{\thesection.}{1em}{} 
\titleformat{\subsection}[block]{\large\bfseries\scshape\centering}{\thesubsection.}{1em}{}
\titleformat{\subsubsection}[block]{\large\bfseries\scshape\centering}{\thesubsubsection.}{1em}{} 

\title{\huge\bfseries Fractional Newton-Raphson Method}

\author[,a,$\star$]{A. Torres-Hernandez  \footnote{E-mail: anthony.torres@ciencias.unam.mx;  ORCID: 0000-0001-6496-9505}}
\affil[a]{\normalsize Department of Physics, Faculty of Science - UNAM, Mexico}

\author[,b,$\star$]{F. Brambila-Paz \footnote{E-mail: fernandobrambila@gmail.com; ORCID: 0000-0001-7896-6460}}
\affil[b]{\normalsize Department of Mathematics, Faculty of Science - UNAM, Mexico}

\date{}

\begin{document}

\maketitle

\thispagestyle{firstpage}

\maketitle\unmarkedfntext{$\star$ Corresponding authors}

\begin{abstract}

The Newton-Raphson (N-R) method is useful to find the roots of a polynomial of degree $ n$, with $ n \in \nset {N}$. However, this method is limited since it diverges for the case in which polynomials only have complex roots if a real initial condition is taken. In the present work, we explain an iterative method that is created using the fractional calculus, which we will call the Fractional Newton-Raphson (F N-R) Method, which has the ability to enter the space of complex numbers given a real initial condition, which allows us to find both the real and complex roots of a polynomial unlike the classical Newton-Raphson method.

\end{abstract}

\section{Newton-Raphson Method}

Let  $\Phi:\nset{R}^n \to \nset{R}^n$ be a function. It is possible to build a sequence $\set{x_i}_{i=0}^\infty$  by defining the following iterative method

\begin{eqnarray}\label{eq:S1-001}
x_{i+1}:=\Phi(x_i),
\end{eqnarray}

if it fulfills that $x_i\to \xi\in \nset{R}^n$ and if the function $\Phi$ is continuous around $\xi$, we obtain that

\begin{eqnarray}\label{eq:S1-002}
\xi=\lim_{i\to \infty}x_{i+1}=\lim_{i\to \infty}\Phi(x_i)=\Phi\left(\lim_{i\to \infty}x_i \right)=\Phi(\xi),
\end{eqnarray}

the above result is the reason by which the method \eqref{eq:S1-001} is known as the \textbf{fixed point method}. Moreover, the function $\Phi$ is called an \textbf{iteration function}.  To understand the nature of the convergence of the iteration function $\Phi$, the following definition is necessary \cite{plato2003concise}:

\begin{definition}
Let $\Phi:\nset{R}^n \to \nset{R}^n$  be an iteration function. The method \eqref{eq:S1-001} for determining $\xi\in \nset{R}^n$ is called \textbf{(locally) convergent}, if there exists $\delta>0$ such that for every initial value

\begin{eqnarray*}
x_0\in B(\xi;\delta):=\set{y\in \nset{R}^n \ : \  \norm{y-\xi}<\delta},
\end{eqnarray*}

it holds that

\begin{eqnarray}\label{eq:S1-003}
\lim_{i \to \infty}\norm{x_i-\xi}\to 0 & \Rightarrow & \lim_{i\to \infty}x_i=\xi.
\end{eqnarray}

\end{definition}

For the one-dimensional case, the N-R method is one of the most used method to find the roots $\xi$ of a function $f: \Omega \subset \nset{R} \to \nset{R} $, that is, $\set{\xi\in \Omega \ : \ f(\xi)=0}$, due to its easy implementation and rapid convergence, the N-R method is expressed in terms of an iteration function $\Phi: \nset{R}\to \nset{R}$, as follows \cite{plato2003concise}:

\begin{eqnarray}\label{eq:S1-004}
x_{i+1}:= \Phi(x_i)=x_i - \left( f^{(1)}(x_i)\right)^{-1} f(x_i),& i=0,1,2,\cdots.
\end{eqnarray}

The N-R method is based on creating a sequence $ \set{x_i}_{i=0}^\infty $ by means of the intersection of the tangent line of the function $ f(x)$ at the $ x_i $ point with the $ x $ axis, if the initial condition $ x_0 $ is close enough to the root $\xi$ then the sequence $ \set{x_i}_{i=0}^\infty $ should be convergent to the root $\xi$ \cite{stoer2013}.  

\begin{figure}[!ht]
\centering
\includegraphics[width=0.6\textwidth, height=0.3\textwidth]{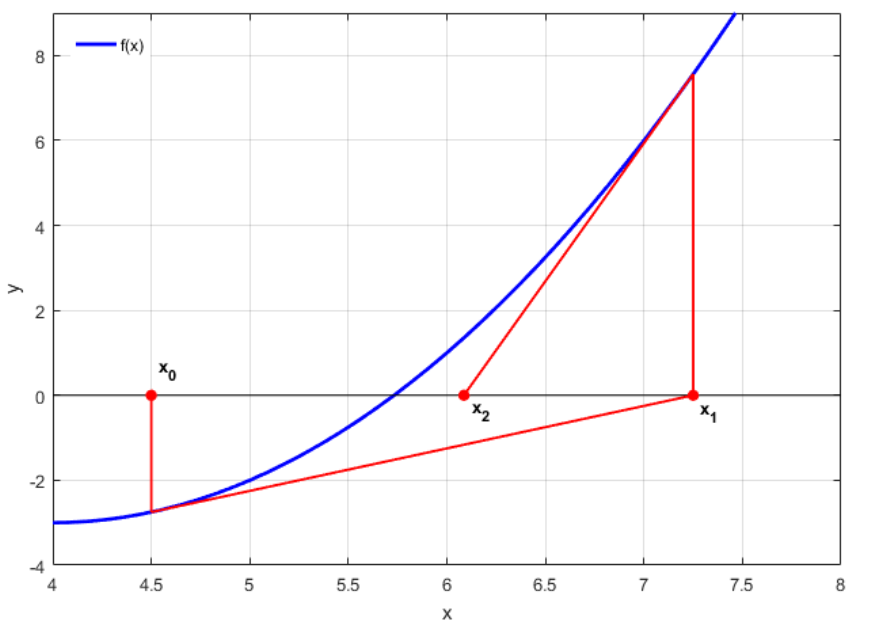}
\caption{ Illustration of the Newton-Raphson method.}
\end{figure}

Before continuing it is necessary to consider the following definition \cite{plato2003concise}:

\begin{definition}
Let $ \Phi: \Omega \subset \nset{R} \to \nset{R} $ be an iteration function with a fixed point $ \xi \in \Omega $. Then the method \eqref{eq:S1-001} is called  \textbf{(locally) convergent of (at least) order $ \boldsymbol{p} $} ($ p \geq 1 $), if there are exists $ \delta> 0 $  and $ C $ a non-negative constant,  with $ C <1 $ if $ p = 1 $, such that for any initial value $ x_0 \in B (\xi; \delta) $ it holds that

\begin{eqnarray}\label{eq:S1-005}
\abs{x_{i+1}-\xi}\leq C \abs{x_i-\xi}^p, & i=0,1,2,\cdots,
\end{eqnarray}

where $ C $ is called convergence factor.

\end{definition}

The order of convergence is usually related to the speed at which the sequence generated by \eqref{eq:S1-001} converges. For the particular case $ p = 1 $ it is said that the method \eqref{eq:S1-001} has an \textbf{order of convergence (at least) linear}, and for the case $p=2$  it is said that the method \eqref{eq:S1-001} has an \textbf{order of convergence (at least) quadratic}. The following theorem, allows characterizing the order of convergence of an iteration function $ \Phi $ with its derivatives \cite{plato2003concise,stoer2013} :

\begin{theorem}\label{teo:c2.01}
Let $ \Phi: \Omega \subset \nset{R} \to \nset {R} $ be an iteration function with a fixed point $ \xi \in \Omega $. Assuming that $\Phi $ is $ p$-times differentiable in $ \xi $ for some $ p \in \nset{N} $, and moreover

\begin{eqnarray}\label{eq:S1-006}
\left\{
\begin{array}{cc}
\ds  \abs{\Phi^{(k)}(\xi)}=0, \ \forall k<p, & \mbox{if }p\geq 2 \vspace{0.1cm}\\
\ds \abs{\Phi^{(1)}(\xi)}<1, & \mbox{if }p=1
\end{array}\right.,
\end{eqnarray}

then $ \Phi $ is (locally) convergent of (at least) order $ p $.

\begin{proof}

Let $\Phi:\nset{R} \to \nset{R}$ be an iteration function, and using the Taylor series expansion of $\Phi$, we obtain two cases:

\begin{itemize}

\item[i)]  Case $p\geq 2:$

\begin{align*}
\Phi(x_i)
= \ds  \Phi(\xi)+  \sum_{k =1}^p \dfrac{\Phi^{(k)}(\xi)}{k !} (x_i-\xi)^k   + o\left( (x_i-\xi)^p\right),
\end{align*}

then

\begin{align*}
\abs{\Phi(x_i)-  \Phi(\xi)}\leq  \sum_{k =1}^p \dfrac{\abs{\Phi^{(k)}(\xi)}}{k !} \abs{x_i-\xi}^k   + o\left( \abs{x_i-\xi}^p\right),
\end{align*}

assuming that $\xi$ is a fixed point of $\Phi$ and that $\abs{\Phi^{(k)}(\xi)}=0 \  \forall k<p$ is fulfilled, the previous expression implies that

\begin{eqnarray*}
\dfrac{\abs{\Phi(x_i)-  \Phi(\xi)}}{\abs{x_i-\xi}^p }=\dfrac{\abs{x_{i+1}- \xi}}{\abs{x_i-\xi}^p }\leq   \dfrac{\abs{\Phi^{(p)}(\xi)}}{p !}   + \dfrac{ o\left( \abs{x_i-\xi}^p\right)}{\abs{x_i-\xi}^p },
\end{eqnarray*}

therefore

\begin{eqnarray*}
\lim_{i\to \infty} \dfrac{\abs{x_{i+1}- \xi}}{\abs{x_i-\xi}^p }\leq   \dfrac{\abs{\Phi^{(p)}(\xi)}}{p !}  ,
\end{eqnarray*}

as a consequence, if the sequence $\set{x_i}_{i=0}^\infty$ generated by \eqref{eq:S1-001} converges to $\xi$, there exists a value $k>0$ such that

\begin{eqnarray*}
\abs{x_{i+1}- \xi}\leq   \dfrac{\abs{\Phi^{(p)}(\xi)}}{p !}\abs{x_i-\xi}^p , & \forall i\geq k,
\end{eqnarray*}

then $ \Phi $ is (locally) convergent of (at least) order $ p $.

\item[ii)] Case $p=1:$

\begin{align*}
\Phi(x_i)
= \Phi(\xi)+ \Phi^{(1)}(\xi) (x_i-\xi) +  o\left( (x_i-\xi)\right)  ,
\end{align*}

then

\begin{align*}
\abs{\Phi(x_i)-\Phi(\xi)}\leq  \abs{\Phi^{(1)}(\xi)} \abs{x_i-\xi }+   o\left( \abs{x_i-\xi} \right) ,
\end{align*}

assuming that $\xi$ is a fixed point of $\Phi$, the previous expression implies that

\begin{eqnarray*}
\dfrac{\abs{\Phi(x_i)-\Phi(\xi)}}{\abs{x_i-\xi}}=\dfrac{\abs{x_{i+1}-\xi}}{\abs{x_i-\xi}}\leq \abs{\Phi^{(1)}(\xi)} +\dfrac{o\left(\abs{x_i-\xi} \right)}{\abs{x_i-\xi}},
\end{eqnarray*}

therefore

\begin{eqnarray*}
\lim_{i\to \infty} \dfrac{\abs{x_{i+1}-\xi}}{\abs{x_i-\xi}}\leq\abs{\Phi^{(1)}(\xi)},
\end{eqnarray*}

as a consequence, if the sequence $\set{x_i}_{i=0}^\infty$ generated by \eqref{eq:S1-001} converges to $\xi$, there exists a value $k>0$ such that

\begin{eqnarray*}
\abs{x_{i+1}-\xi}\leq \abs{\Phi^{(1)}(\xi)}\abs{x_i-\xi}, & \forall i\geq k,
\end{eqnarray*}

then considering $m\geq 1$

\begin{align*}
\abs{x_{i+m}-\xi}\leq &\abs{\Phi^{(1)}(\xi)}\abs{x_{i+m-1}-\xi}\leq \abs{\Phi^{(1)}(\xi)}^2\abs{x_{i+m-2}-\xi} \leq \cdots \leq\abs{\Phi^{(1)}(\xi)}^{m}\abs{x_{i}-\xi},
\end{align*}

and assuming that $\abs{\Phi^{(1)}(\xi)}<1$ is fulfilled

\begin{eqnarray*}
\lim_{m \to \infty}\abs{x_{i+m}-\xi}\leq \lim_{m \to \infty}\abs{\Phi^{(1)}(\xi)}^{m}\abs{x_{i}-\xi} \to 0,
\end{eqnarray*}

then $ \Phi $ is (locally) convergent of order (at least) linear.

\end{itemize}

\end{proof}

\end{theorem}

The N-R method is characterized by having an order of convergence at least quadratic for the case where $ f^{(1)}(\xi) \neq 0$, but if to the previous case it is added that $f^{(2)}(\xi)=0$, then the N-R method presents an order of convergence at least cubic. On other hand, for the case where the function $f$ has a root $\xi$ with a certain algebraic multiplicity $ m\geq 2$, that is, 

\begin{eqnarray*}
\begin{array}{cc}
f(x)=(x-\xi)^mg(x), & g(\xi)\neq 0,
\end{array}
\end{eqnarray*}

the N-R method presents an order of convergence at least linear \cite{plato2003concise}. The aforementioned may be formalized by the following proposition:

\begin{proposition}\label{prop:01}
Let $f:\Omega \subset \nset{R} \to \nset{R}$ be a function with a zero $\xi \in \Omega$. Then the iteration function $\Phi$ of the N-R method, given by \eqref{eq:S1-004}, fulfills the following condition:

\begin{eqnarray}\label{eq:S1-007}
\abs{x_{i+1}- \xi}\leq   \dfrac{\abs{\Phi^{(p)}(\xi)}}{p !}\abs{x_i-\xi}^p,
\end{eqnarray}

where

\begin{eqnarray}
p=\left\{
\begin{array}{cc}
1 ,& \mbox{ if } f(x)=(x-\xi)^mg(\xi) \mbox{ with } g(\xi)\neq 0 \mbox{ and }m\geq 2 \vspace{0.1cm}\\
2, & \mbox{ if } f^{(1)}(\xi)\neq 0 \vspace{0.1cm}\\
3,& \mbox{ if } f^{(1)}(\xi)\neq 0 \mbox{ and } f^{(2)}(\xi)=0
\end{array}\right. .
\end{eqnarray}

\begin{proof}
Considering that the form of the function $f$ is not explicitly determined, it is possible to consider two possibilities:

\begin{itemize}
\item[i)] Assuming the function may be written as $f(x)=(x-\xi)^mg(\xi)$  with  $g(\xi)\neq 0$  and $m\geq 2$, then

\begin{eqnarray*}
f^{(1)}(x)=(x-\xi)^{m-1}\left[(x-\xi) g^{(1)}(x)+mg(x)\right],
\end{eqnarray*}

as a consequence, the iteration function of N-R method takes the following form

\begin{eqnarray*}
\Phi(x)= x-(x-\xi)h(x)g(x),
\end{eqnarray*}

with

\begin{eqnarray*}
h(x)=\left[(x-\xi) g^{(1)}(x)+mg(x)\right]^{-1},
\end{eqnarray*}

then

\begin{align*}
\Phi^{(1)}(x)= 1- h(x)\left[(x-\xi)g^{(1)}(x)+g(x)  \right]-(x-\xi)h^{(1)}(x)g(x),
\end{align*}

where

\begin{eqnarray*}
h^{(1)}(x)=-\left[(x-\xi) g^{(1)}(x)+mg(x)\right]^{-2}\left[\left( 1+m \right) g^{(1)}(x)+(x-\xi) g^{(2)}(x) \right],
\end{eqnarray*}

therefore

\begin{align}
\lim_{x\to \xi}\abs{\Phi^{(1)}(x)}= \abs{1-h(\xi)g(\xi)}=\abs{1-\dfrac{1}{m}}<1,
\end{align}

and from the \textbf{Theorem \ref{teo:c2.01}}, the N-R method has an order of convergence at least linear, that is, the N-R method fulfills the equation \eqref{eq:S1-007} with $p=1$.

\item[ii)] Assuming that $f(x)\neq (x-\xi)^mg(\xi)$  with  $g(\xi)\neq 0$  and $m\geq 2$, the first derivative of the iteration function of N-R method takes the following form

\begin{align*}
\Phi^{(1)}(x)=\left(f^{(1)}(x) \right)^{-2}f(x)f^{(2)}(x),
\end{align*}

and if it fulfills that $f^{(1)}(\xi)\neq 0$, then

\begin{align}
\lim_{x\to \xi}\abs{\Phi^{(1)}(x)}= 0,
\end{align}

and from the \textbf{Theorem \ref{teo:c2.01}}, the N-R method has an order of convergence at least quadratic, that is, the N-R method fulfills the equation \eqref{eq:S1-007} with $p=2$. On other hand, the second derivative of the iteration function of N-R method takes the following form

\begin{align*}
\Phi^{(2)}(x)= \left(f^{(1)}(x) \right)^{-1}f^{(2)}(x)+ f(x)\left[\left(f^{(1)}(x) \right)^{-2} f^{(3)}(x)-2 \left(f^{(1)}(x) \right)^{-3}\left( f^{(2)}(x) \right)^2 \right],
\end{align*}

and if it fulfills that $f^{(1)}(\xi)\neq 0$ and $f^{(2)}(\xi)=0$, then

\begin{align}
\lim_{x\to \xi}\abs{\Phi^{(1)}(x)}=\lim_{x\to \xi}\abs{\Phi^{(2)}(x)}= 0,
\end{align}

and from the \textbf{Theorem \ref{teo:c2.01}}, the N-R method has an order of convergence at least cubic, that is, the N-R method fulfills the equation \eqref{eq:S1-007} with $p=3$.

\end{itemize}

\end{proof}

\end{proposition}

The previous proposition, illustrates two important points that are worth mentioning when using the N-R method to find the zeros of a function $f$:

\begin{itemize}
\item[i)] When it is not evident, unless it is explicitly specified that the function $f$ has no roots of algebraic multiplicity $m\geq 2$, technically there exists the possibility that the N-R method has an order of convergence at least linear, that is, the N-R method may fulfill the equation \eqref{eq:S1-007} with $p\geq 1$.

\item[ii)] Due that the N-R method is a local iterative method, even if it proves that for a root $\xi\in \Omega$ the method has an order of convergence at least linear, this does not rule out that for the same function  $f$ it may present a higher order of convergence over the same region $\Omega$. As an example of the above, we may consider the following function

\begin{eqnarray*}
f(x)=(x-\eta)(x-\xi)^mg(x), & g(\eta)\neq g(\xi)\neq 0,
\end{eqnarray*}

with $ \eta,\xi \in \Omega$, $\abs{\eta-\xi}<\epsilon$, and  $m\geq 2$.

\end{itemize}

The previous points are important, because when the N-R method is implemented in a function $f$, the zeros of the function are assumed to be unknown, and their algebraic multiplicities $m\geq 2$, in case they exist, are also unknown. With the above in mind, the following corollary is obtained, which is derived from the \textbf{Theorem \ref{teo:c2.01}}:

\begin{corollary}\label{cor:2-001}
Let $\Phi:\nset{R} \to \nset{R}$ be an iteration function. If $\Phi$ defines a sequence $\set{x_i}_{i=0}^\infty$ such that $x_i\to \xi$, and if the following condition is fulfilled

\begin{eqnarray}\label{eq:S1-008}
\lim_{x\to \xi}\abs{\Phi^{(1)}(x)}\neq 0,
\end{eqnarray}

then $\Phi$ has an order of convergence (at least) linear  in $B(\xi;\delta)$.
\end{corollary}

\section{Fractional Calculus}

The fractional calculus is a mathematical analysis branch whose applications have been increasing since the end of the XX century and beginnings of the XXI century \cite{brambila2017fractal,martinez2017applications2}, the fractional calculus arises around 1695 due to Leibniz's notation for the derivatives of integer order

\begin{eqnarray*}
f^{(n)}(x):=\dfrac{d^n}{dx^n}f(x), & n\in \nset{N},
\end{eqnarray*}

thanks to this notation L'Hopital could ask in a letter to Leibniz about the interpretation of taking $ n = 1/2 $ in a derivative, since at that moment Leibniz could not give a physical or geometrical interpretation to this question, he simply answered L'Hopital in a letter, \com{$\dots$ is an apparent paradox of which, one day, useful consequences will be drawn} \cite{miller93}. The name of fractional calculus comes from a historical question since in this branch of mathematical analysis it is studied the derivatives and integrals of a certain order $\alpha$, with $\alpha \in \nset{R}$ or $\nset{C}$. 

Currently, the fractional calculus does not have a unified definition of what is considered a fractional derivative, because one of the conditions required to consider an expression as a fractional derivative is to recover the results of conventional calculus when the order $\alpha \to n$, with $n \in \nset{N}$ \cite{oldham74}, among the most common definitions of fractional derivatives are the  Riemann-Liouville (R-L) fractional derivative and the Caputo fractional derivative \cite{hilfer00}, the latter is usually the most studied since the Caputo fractional derivative allows us a physical interpretation to problems with initial conditions, this derivative fulfills the property of the classical calculus that the derivative of a constant is null regardless of the order $\alpha$ of the derivative, however this does not occur with the R-L fractional derivative, and this characteristic can be used to solve nonlinear systems \cite{torres2020reduction,torres2020nonlinear}.

Unlike the Caputo fractional derivative, the R-L fractional derivative does not allow for a physical interpretation to the problems with initial condition because its use induces  fractional initial conditions, however the fact that this derivative does not cancel the constants for $\alpha$, with $\alpha\notin \nset{N}$, allows to obtain a \com{spectrum} of the behavior of the constants for different orders of the derivative, which is not possible with conventional calculus.

\subsection{Introduction to the Riemann-Liouville Fractional Derivative }

One of the key pieces in the study of fractional calculus is the iterated integral, which is defined as follows \cite{hilfer00}

\begin{definition}
Let $ L_{loc} ^ 1 (a, b) $, the space of locally integrable functions in the interval $ (a, b) $. If $ f $ is a function such that $ f \in L_ {loc} ^ 1 (a, \infty) $, then the $n$-th iterated integral of the function $ f $ is given by 

\begin{eqnarray}\label{eq:S2-001}
\begin{array}{c}
\ds \ifr{}{a}{I}{x}{n} f(x)=\ifr{}{a}{I}{x}{}\left(\ifr{}{a}{I}{x}{n-1} f(x)  \right)=\frac{1}{(n-1)!}\int_a^x(x-t)^{n-1}f(t)dt,
\end{array}
\end{eqnarray}

where

\begin{eqnarray*}
\ifr{}{a}{I}{x}{} f(x):=\int_a^x f(t)dt.
\end{eqnarray*}

\end{definition}

Considerate that $ (n-1)! = \gam{n} $
, a generalization of \eqref{eq:S2-001} may be obtained for an arbitrary order $ \alpha> 0 $

\begin{eqnarray}\label{eq:S2-002}
\ifr{}{a}{I}{x}{\alpha} f(x)=\dfrac{1}{\gam{\alpha}}\int_a^x(x-t)^{\alpha-1}f(t)dt,
\end{eqnarray}

similarly, if $ f \in L_{loc} ^ 1 (- \infty, b) $, we may define

\begin{eqnarray}\label{eq:S2-003}
\ifr{}{x}{I}{b}{\alpha} f(x)=\dfrac{1}{\gam{\alpha}}\int_x^b(t-x)^{\alpha -1}f(t)dt,
\end{eqnarray} 

the equations \eqref{eq:S2-002} and \eqref{eq:S2-003} correspond to the definitions of \textbf{right and left Riemann-Liouville fractional integral}, respectively. The fractional integrals fulfill the  \textbf{semigroup property}, which is given in the following proposition \cite{hilfer00}:

\begin{proposition}
Let $ f $ be a function. If $ f \in L_{loc} ^ 1 (a, \infty) $, then the fractional integrals of $ f $ fulfill that

\begin{eqnarray}\label{eq:S2-004}
\ifr{}{a}{I}{x}{\alpha} \ifr{}{a}{I}{x}{\beta}f(x) = \ifr{}{a}{I}{x}{\alpha + \beta}f(x),& \alpha,\beta>0.
\end{eqnarray}

\end{proposition}

From the previous result, and considering that the operator $ d / dx $  is the inverse operator to the left of the operator $ \ifr {}{a}{I}{x}{} $, any integral $ \alpha$-th of a function $ f \in L_{loc} ^ 1 (a, \infty) $ may be written as

\begin{eqnarray}\label{eq:S2-005}
\ifr{}{a}{I}{x}{\alpha}f(x)=\dfrac{d^n}{dx^n}\left( \ifr{}{a}{I}{x}{n}\ifr{}{a}{I}{x}{\alpha}f(x) \right)=\dfrac{d^n}{dx^n}\left( \ifr{}{a}{I}{x}{n+\alpha}f(x)\right).
\end{eqnarray}
 
Considering \eqref{eq:S2-002} and \eqref{eq:S2-005}, we can built the \textbf{Riemann-Liouville fractional derivative} $\ifr{}{a}{D}{x}{\alpha}$, as follows \cite{hilfer00,kilbas2006theory}:

\begin{eqnarray}\label{eq:S2-006}
\normalsize
\begin{array}{c}
\ifr{}{a}{D}{x}{\alpha}f(x) := \left\{
\begin{array}{cc}
\ds \ifr{}{a}{I}{x}{-\alpha}f(x), &\mbox{if }\alpha<0\\  
\ds \dfrac{d^n}{dx^n}\left( \ifr{}{a}{I}{x}{n-\alpha}f(x)\right), & \mbox{if }\alpha\geq 0
\end{array}
\right. ,
\end{array}
\end{eqnarray}

where  $ n = \lceil \alpha \rceil$, then applying the  operator \eqref{eq:S2-006} to the  function $ x^{\mu} $, with  $\alpha \in \nset{R}\setminus\nset {Z} $ and $\mu\geq 0$, we obtain the following result

\begin{eqnarray}\label{eq:S2-007}
\ifr{}{0}{D}{x}{\alpha}x^\mu = 
 \dfrac{\gam{\mu+1}}{\gam{\mu-\alpha+1}}x^{\mu-\alpha}.
\end{eqnarray}

\subsection{Introduction to the Caputo Fractional Derivative}

Michele Caputo (1969) published a book and introduced a new definition of fractional derivative, he created this definition with the objective of modeling anomalous diffusion phenomena. The definition of Caputo had already been discovered independently by Gerasimov (1948). This fractional derivative is of the utmost importance since it allows us to give a physical interpretation of the initial value problems, moreover to being used to model fractional time. In some texts, it is known as the fractional derivative of Gerasimov-Caputo.

Let $ f $ be a function, such that $ f $ is $ n$-times differentiable with $ f ^{(n)} \in L_{loc}^ 1 (a, b) $, then the \textbf{(right) Caputo fractional derivative} is defined as \cite{kilbas2006theory}

\begin{align}\label{eq:S2-008}
\ifr{C}{a}{D}{x}{\alpha}f(x):= &\ifr{}{a}{I}{x}{n-\alpha}\left( \der{d}{x}{n} f(x)\right) = \dfrac{1}{\gam{n-\alpha}}\int_{a}^{x} (x-t)^{n-\alpha -1} f^{(n)}(t)dt , 
\end{align}

where $n=\lceil \alpha \rceil$. It should be mentioned that the Caputo fractional derivative behaves as the inverse operator to the left of the Riemann-Liouville fractional integral, that is,

\begin{eqnarray*}
\ifr{C}{a}{D}{x}{\alpha}(\ifr{}{a}{I}{x}{\alpha}f(x))=f(x).
\end{eqnarray*}

On the other hand, the relation between the fractional derivatives of Caputo and Riemann-Liouville is given by the following expression \cite{kilbas2006theory}

\begin{eqnarray*}
\ifr{C}{a}{D}{x}{\alpha}f(x)=\ifr{}{a}{D}{x}{\alpha}\left(f(x)-\sum_{k=0}^{n-1}\dfrac{f^{(k)}(a)}{k!}(x-a)^k\right), 
\end{eqnarray*}

then, if $f^{(k)}(a)=0 \ \ \forall k<n$, we obtain

\begin{eqnarray*}
\ifr{C}{a}{D}{x}{\alpha}f(x)=\ifr{}{a}{D}{x}{\alpha}f(x), 
\end{eqnarray*}

considering the previous particular case, it is possible to unify the definitions of R-L fractional integral and Caputo fractional derivative as follows

\begin{eqnarray}\label{eq:S2-009}
\begin{array}{c}
\ifr{C}{a}{D}{x}{\alpha}f(x) := \left\{
\begin{array}{cc}
\ds \ifr{}{a}{I}{x}{-\alpha}f(x), &\mbox{if }\alpha<0\\  
\ds \ifr{}{a}{I}{x}{n-\alpha}\left( \der{d}{x}{n} f(x)\right) , & \mbox{if }\alpha\geq 0
\end{array}
\right. .
\end{array}
\end{eqnarray}

\section{Fractional Newton-Raphson Method}

Let $\nset{P}_n(\nset{R})$ be the space of polynomials of degree less than or equal to $ n\in \nset{N} $ with real coefficients. The N-R method is useful for finding the roots of a function $ f$. However, this method is limited because it cannot find roots $ \xi \in \nset{C} \setminus \nset {R} $, if the sequence $ \set{x_i}_{i = 0} ^ \infty $ generated by \eqref{eq:S1-004} has an initial condition $ x_0 \in \nset{R} $. To solve this problem and develop a method that has the ability to find roots, both real and complex, of a polynomial if the initial condition $ x_0 $ is real, we propose a new method, which consists of Newton-Raphson method with the implementation of the fractional derivative. Before continuing, it is necessary to define the following notation

\begin{eqnarray}\label{eq:S3-001}
f^{(\alpha)}(x):=\dfrac{d^{\alpha}}{dx^\alpha}f(x) ,
\end{eqnarray}

where the operator $ d^ \alpha / dx^\alpha $ denotes any fractional derivative, applied on the variable $x$, that fulfills the following condition of continuity respect to the order of the derivative

\begin{eqnarray}\label{eq:S3-002}
\lim_{\alpha\to 1}f^{(\alpha)}(x)=f^{(1)}(x).
\end{eqnarray}

Considering a function $\Phi:(\nset{R}\setminus \nset{Z})\times \nset{C} \to \nset{C}$. Then, using as a basis the idea of the N-R method \eqref{eq:S1-004}, and considering any fractional derivative that fulfills the condition \eqref{eq:S3-002}, we can define the \textbf{Fractional Newton-Raphson method} as follows:

\begin{eqnarray}\label{eq:S3-003}
\begin{array}{cc}
x_{i+1}:= \Phi\left( \alpha, x_i \right)=x_i-\left(f^{\left(\alpha\right)}(x_i) \right)^{-1}  f(x_i),& i=0,1,2,\cdots.
\end{array}
\end{eqnarray}

For the above expression to make sense, due to the part of the integral operator that fractional derivatives usually have, and  that the F N-R method can be used in a wide variety of functions \cite{torreshern2020}, we consider in the expression \eqref{eq:S3-003} that the fractional derivative is obtained for a real variable $x$, and if the result allows it, this variable is subsequently made to tend to a complex variable $x_i$,  that is,

\begin{eqnarray}
f^{\left(\alpha\right)}(x_i):=f^{(\alpha)}(x)\bigg{|}_{ x\longrightarrow x_i}  , & x\in \nset{R}, & x_i\in \nset{C}.
\end{eqnarray}

It should be mentioned that in general, in the F N-R method  $\abs{\Phi^{(1)}(\alpha,\xi)}\neq 0$ if $f(\xi)=0$, and from the \textbf{Corollary \ref{cor:2-001}}, the \textbf{Proposition \ref{prop:01}} and the condition \eqref{eq:S3-002}, any sequence $ \set{x_i} _ {i = 0} ^ \infty $ generated by the iterative method \eqref {eq:S3-003} has an order of convergence at least linear, that is, the F N-R method may fulfill the equation \eqref{eq:S1-007} with $p\geq 1$,  which becomes more evident when considering $\alpha\in [1-\epsilon,1+\epsilon]\setminus \set{1}$.  

To understand why the F N-R method, if $f\in \nset{P}_n(\nset{R}$), has the ability to enter the complex space using a real initial condition unlike the classical N-R method, it is enough to observe the R-L fractional derivative \eqref{eq:S2-007},  with $\alpha=1/2$, of the constant function $f_0(x)=x^0$ and the identity function $f_1(x)=x^1$:

\begin{eqnarray*}
\ifr{}{0}{D}{x}{1/2}f_0(x) = 
\dfrac{\gam{1}}{\gam{1/2}}  x^{-1/2} ,  & \ifr{}{0}{D}{x}{1/2}f_1(x) = 
\dfrac{\gam{2}}{\gam{3/2}}  x^{1/2}.
\end{eqnarray*}

\begin{figure}[!ht]
\centering
\begin{subfigure}[c]{0.49\linewidth}
\includegraphics[width=\linewidth,height=0.42\linewidth]{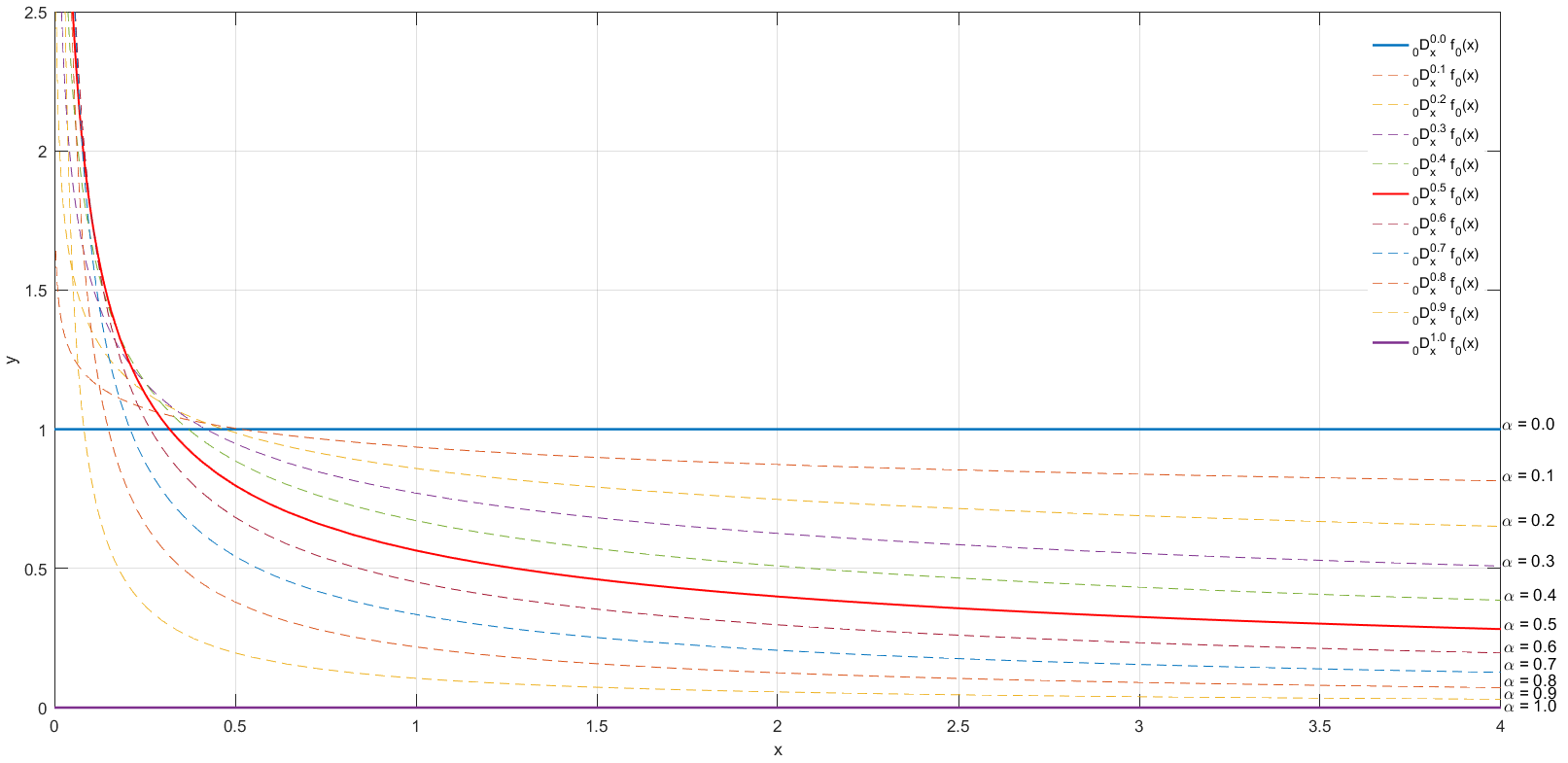}
\end{subfigure}
\begin{subfigure}[c]{0.49\linewidth}
\includegraphics[width=\linewidth,height=0.42\linewidth]{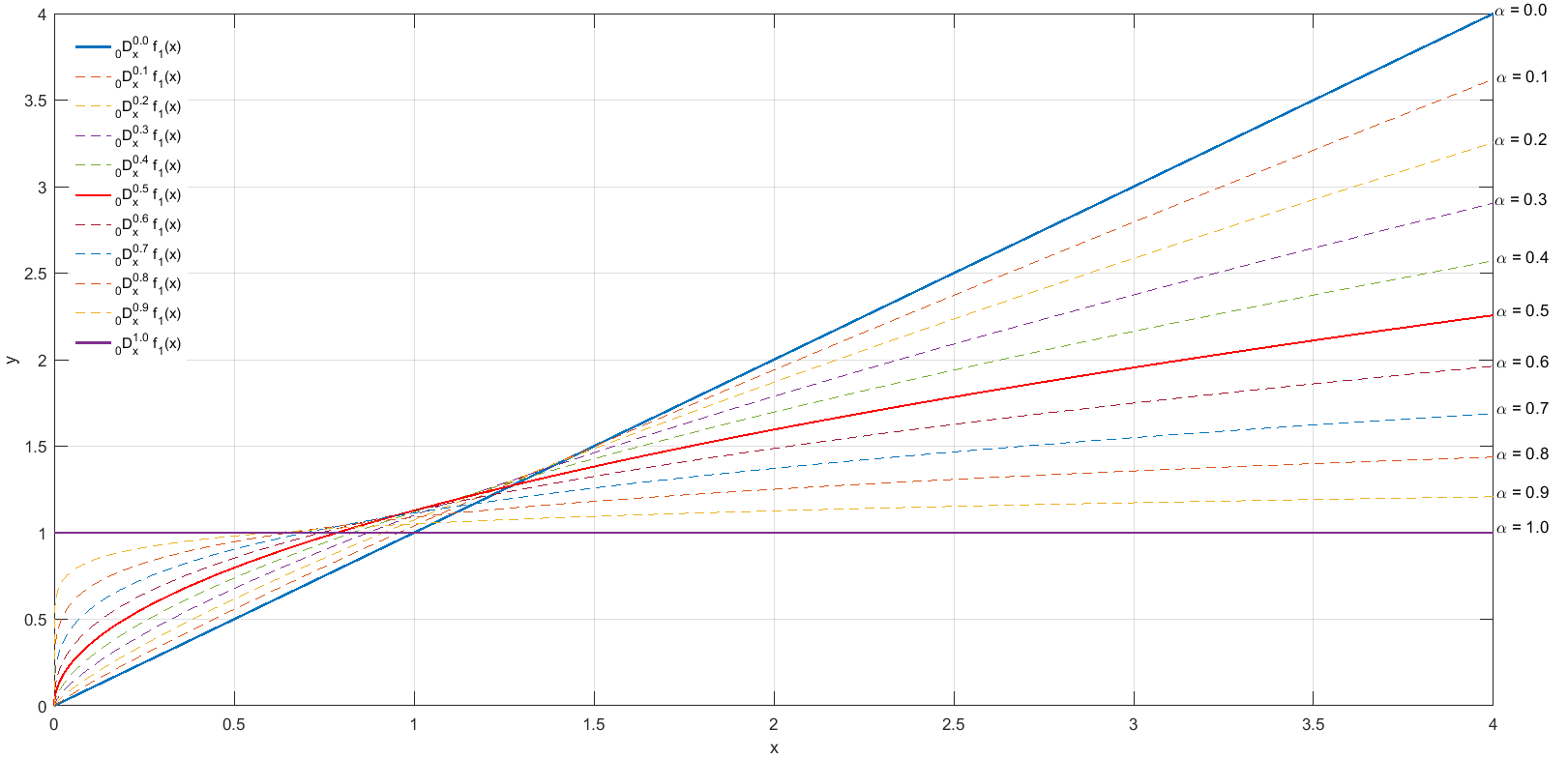}
\end{subfigure}
\caption{The R-L fractional derivatives of $f_0(x)$ and $f_1(x)$, with $\alpha\in[0,1]$.}
\end{figure}

For polynomials of degree $n\geq 1$, in the F N-R method the initial condition $x_0$ must be taken different to zero, as a consequence of the R-L fractional derivative of order $\alpha$, with $\alpha\notin \nset{Z}$, of the constants are proportional to the function $x^{-\alpha}$. When using the F N-R method, with the R-L fractional derivative, on a function $f \in \nset{P}_{n}(\nset{R})$, presents among its behaviors, the following particular cases depending on the initial condition $ x_0 $: 

\begin{itemize}
\item[i)] If we take an initial condition $ x_0> 0 $, the sequence $ \set{x_i}_{i = 0} ^ \infty $ may be divided into three parts, this occurs because it may exists a value $ M \in \nset{N} $  for which $ \set{x_i}_{i = 0}^{M-1} \subset \nset{R}_{>0} $ with $ \set{x_M} \subset \nset{R}_{<0}$, in consequence $ \set {x_i} _ {i \geq M + 1} \subset \nset{C} $.

\item[ii)] On the other hand, if  we take an initial condition  $ x_0 <0 $,  the sequence $ \set {x_i}_ {i = 0} ^ \infty $ may  be divided into two parts, $\set{x_0}\subset \nset{R}_{<0}$ and  $\set{x_i}_{i\geq 1}\subset \nset{C}$.
\end{itemize}

\subsection{Advantages of the Fractional Newton-Raphson Method}

One of the main advantages of the F N-R method is that the initial condition $x_0$ can be left fixed, and so vary the order $\alpha$ of the derivative to obtain both real and complex roots of a polynomial. Due that the  order $\alpha $ of the derivative is varied, different values of $\alpha$ can throw the same root but with a different number of iteration, so to optimize the method, it is possible to implement a filter in which once we have obtained the roots, only those whose orders of the derivatives have generated a smaller number of iterations are extracted.

\begin{figure}[!ht]
\centering
\includegraphics[width=0.6\textwidth, height=0.35\textwidth]{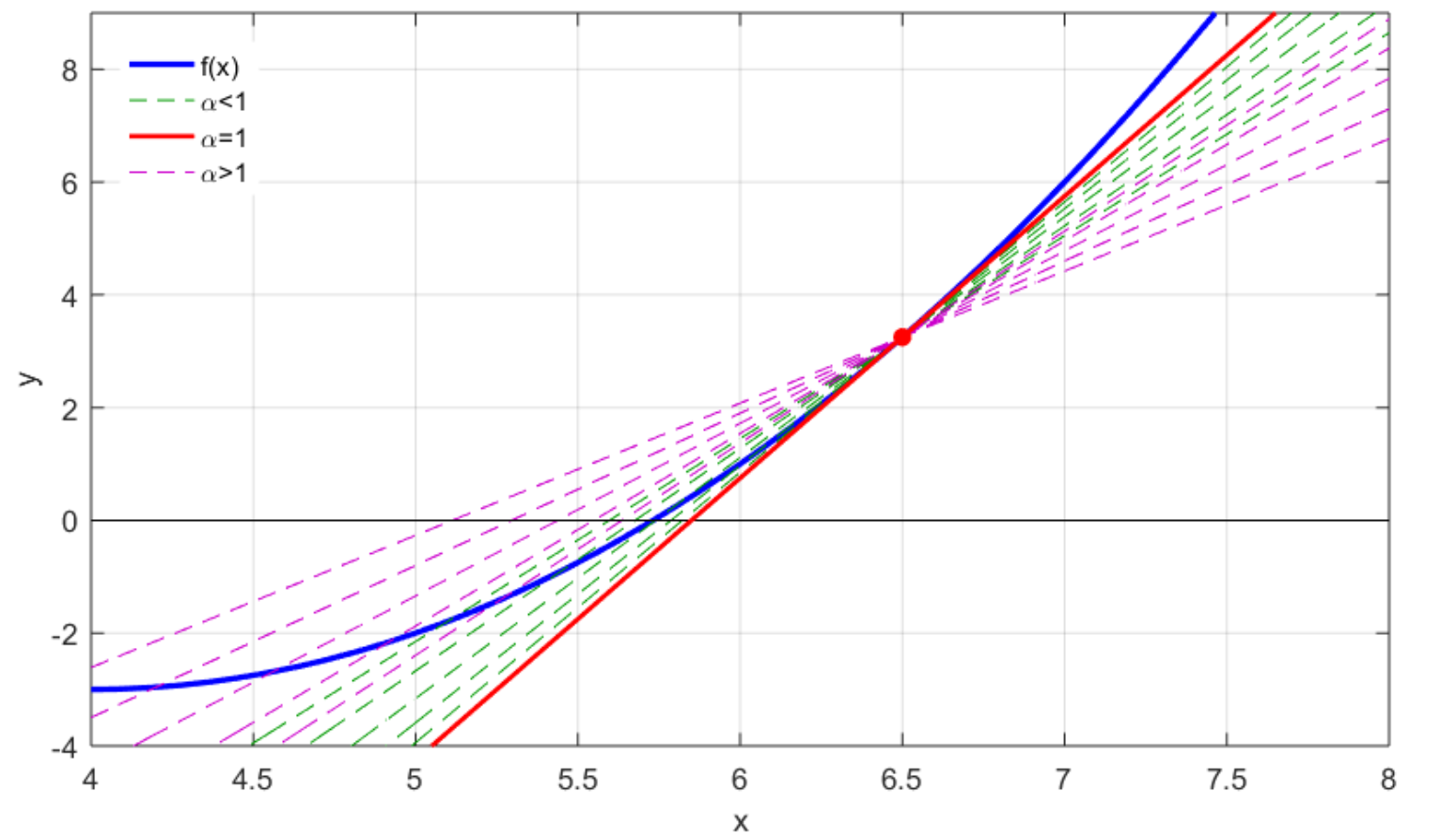}
\caption{Illustration of some lines generated by the F N-R method.}
\end{figure}

\begin{figure}[!ht]
        \begin{subfigure}[c]{0.24\textwidth}
        \centering
 \includegraphics[width=\textwidth, height=1.1\textwidth]{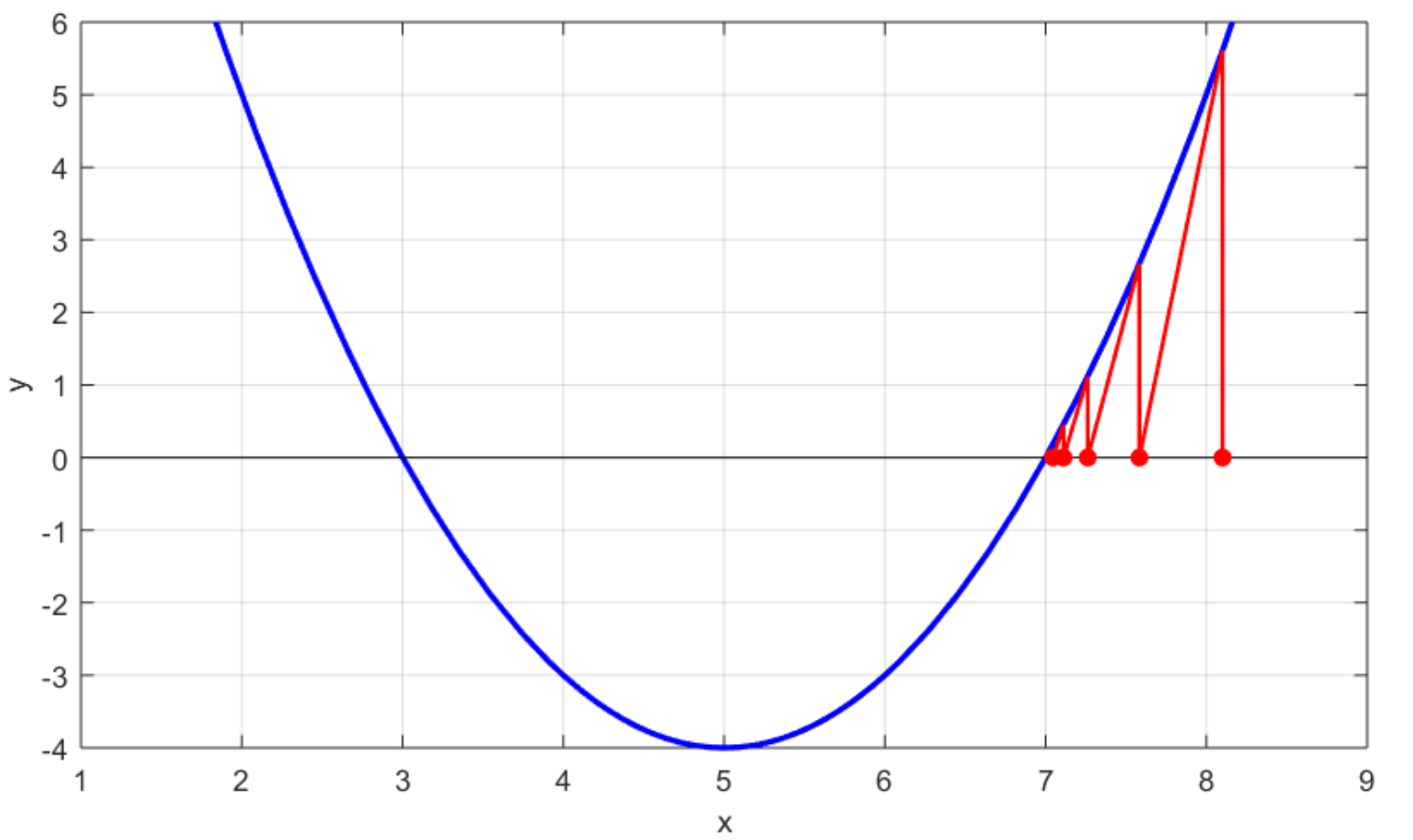}      
    \caption*{a) $\alpha=-0.77$}
    \end{subfigure}
        \begin{subfigure}[c]{0.24\textwidth}
        \centering
 \includegraphics[width=\textwidth, height=1.1\textwidth]{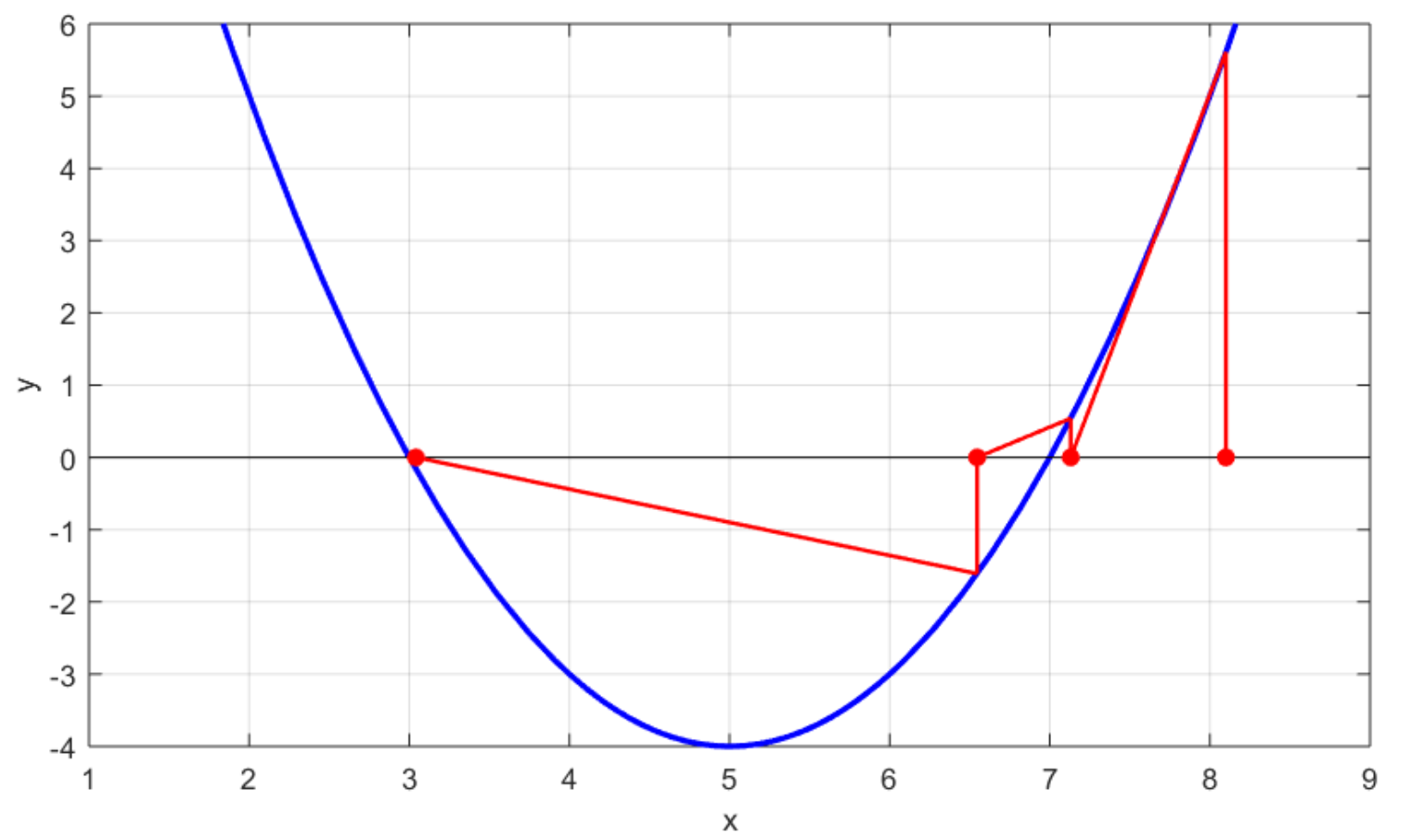}     
    \caption*{b) $\alpha=-0.32$}
    \end{subfigure}    
    \centering
    \begin{subfigure}[c]{0.24\textwidth}
    \centering
 \includegraphics[width=\textwidth, height=1.1\textwidth]{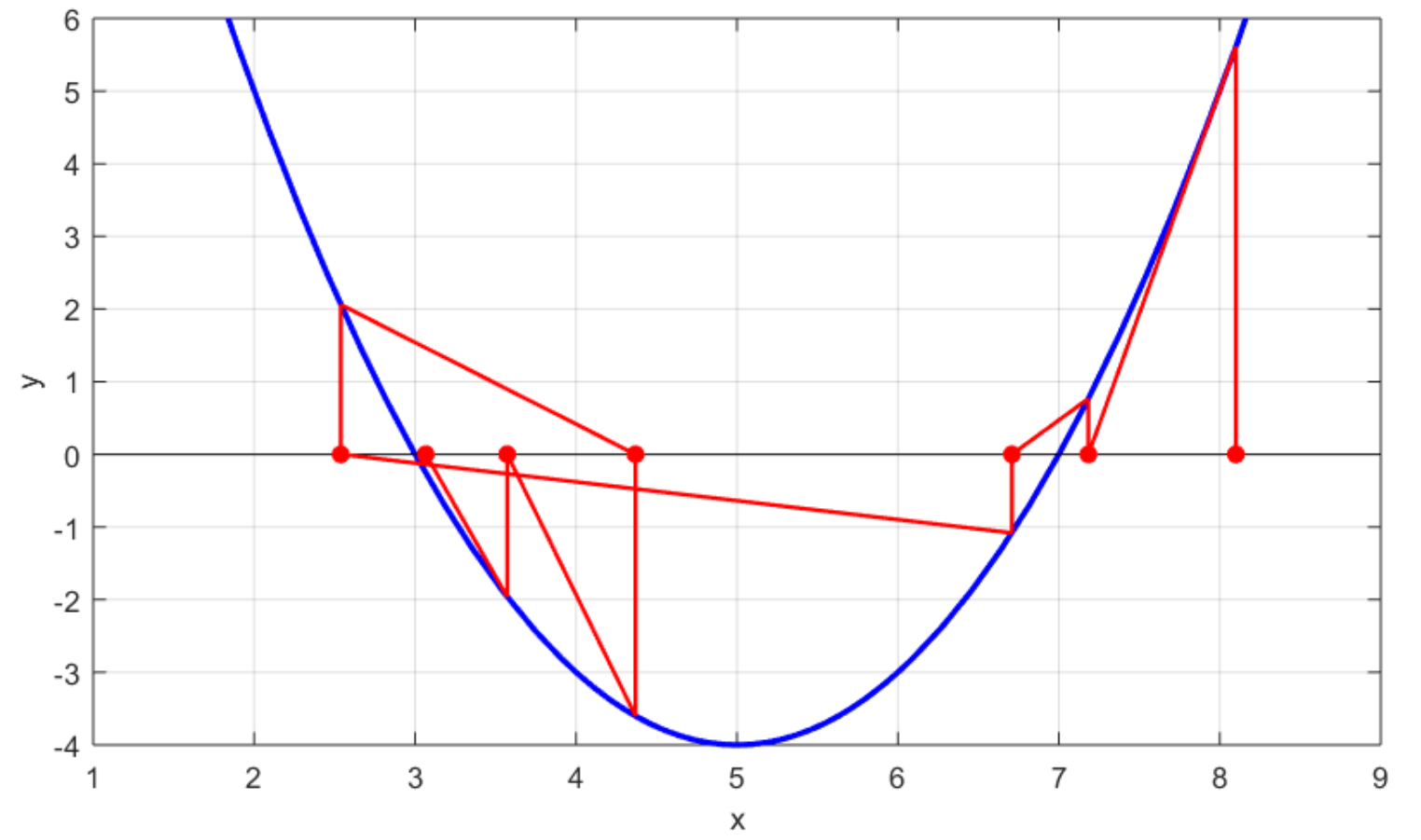}      
    \caption*{c) $\alpha=0.19$}
    \end{subfigure}
    \begin{subfigure}[c]{0.24\textwidth}
    \centering
 \includegraphics[width=\textwidth, height=1.1\textwidth]{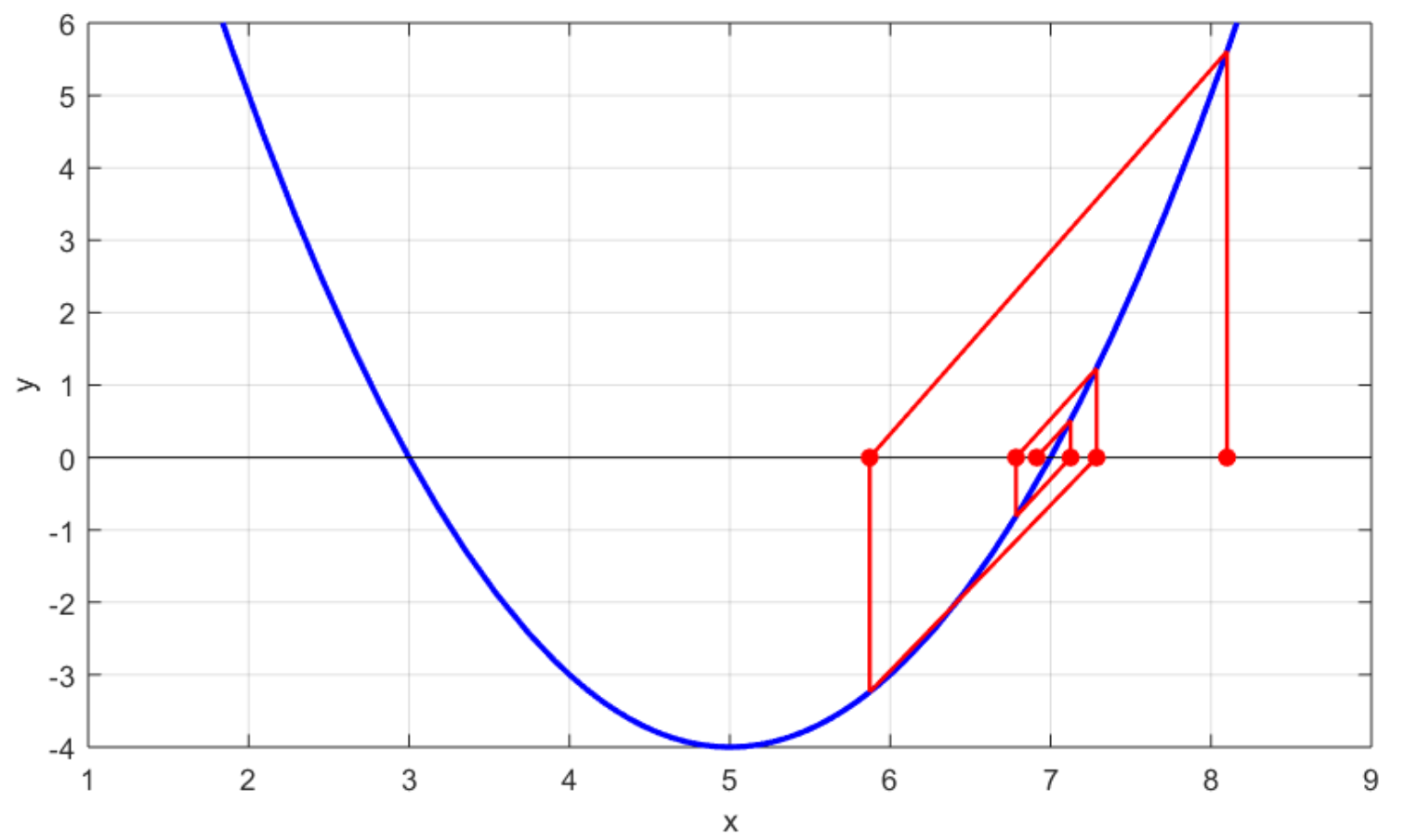}      
    \caption*{d) $\alpha=1.87$}
    \end{subfigure}        
        \caption{llustrations of some trajectories generated by the F N-R method for the same initial condition $ x_0 $ but with different values of $ \alpha $.}\label{fig:02}
\end{figure}

Another advantage is a consequence that the method provides complex roots, so once a root is obtained it is enough to obtain its conjugate complex to obtain another root, in essence, it could be considered that they extract two roots with the same order of the derivative and the same number of iterations. The method does not guarantee that all roots of the polynomial are found by leaving an initial condition fixed and by varying the orders $ \alpha $ of the derivative, as in the classical N-R method, finding the roots will depend on giving an appropriate initial condition.

\subsection{Results of the Fractional Newton-Raphson Method}

The following examples are solved using the R-L fractional derivative \eqref{eq:S2-007}. Instructions for implementing the F N-R method, along with information to provide values $\alpha \in [0.7,1.3]\setminus \set{1}$ are found in the reference \cite{torreshern2020}.  For rounding reasons, for the examples the following function is defined

\begin{eqnarray}\label{eq:S3-004}
\rnd{x}{m}:=\left\{
\begin{array}{cc}
\re{x},& \mbox{ if \hspace{0.1cm}} \abs{\im{x}}\leq 10^{-m}\vspace{0.1cm}\\
x,& \mbox{ if \hspace{0.1cm}} \abs{\im{x}}> 10^{-m}\vspace{0.1cm}\\
\end{array}\right..
\end{eqnarray}

Combining the function \eqref{eq:S3-004} with the method \eqref{eq:S3-003}, the following iterative method is defined

\begin{eqnarray}\label{eq:S3-005}
x_{i+1}:=\rnd{\Phi(\alpha, x_i)}{5}, & i=0,1,2\cdots.
\end{eqnarray}

\begin{example}

Let $f$ be a function, with

\begin{align*}
f(x)=&-64.23x^{14}-72.74x^{13}-61.66x^{12}+32.26x^{11}\\
&+32.3x^{10}-41.37x^9+20.18x^8+4.32x^7\\
&-5.67x^6+17.41x^5-78.6x^4-48.27x^3\\
&
-19.31x^2+77.92x-45.03.
\end{align*}

Then the initial condition $x_0=0.68$ is chosen to use the iterative method  given by \eqref{eq:S3-005}. Consequently, we obtain the results of the Table \ref{tab:01}.

\begin{table}[!ht]
\centering
$
\begin{array}{c|ccccc}
\toprule
&\alpha& x_n&\norm{x_n - x_{n-1} }_2  &\norm{f\left(x_n \right)}_2& n \\ 
\midrule
    1     & 0.97399 & 0.89785306 - 0.29205148i & 1.65245E-06 & 9.65908E-05 & 17 \\
    2     & 0.97455 & 0.53497472 - 0.82703363i & 6.75426E-07 & 7.71022E-05 & 19 \\
    3     & 0.97631 & -0.07482893 - 1.0188353i & 2.60768E-07 & 3.30325E-05 & 17 \\
    4     & 0.98539 & -0.6673652 - 1.16572645i & 1.06301E-07 & 3.96105E-05 & 26 \\
    5     & 0.99275 & -0.67766753 - 0.66659064i & 1.34536E-07 & 6.10746E-06 & 18 \\
    6     & 1.00575 & -0.07482895 + 1.01883532i & 3.83275E-07 & 1.21870E-05 & 22 \\
    7     & 1.00623 & -0.6673652 + 1.16572646i & 1.66433E-07 & 5.94548E-05 & 32 \\
    8     & 1.00635 & -0.67766754 + 0.66659064i & 1.50333E-07 & 2.08269E-06 & 19 \\
    9     & 1.00643 & 0.53497473 + 0.82703358i & 1.70294E-07 & 9.71477E-06 & 21 \\
    10    & 1.03515 & -1.09479584 - 0.25179059i & 1.44222E-07 & 5.21878E-05 & 26 \\
    11    & 1.15163 & -1.09479581 + 0.25179059i & 6.08276E-08 & 9.98791E-05 & 25 \\
    12    & 1.15715 & 0.51558361 - 0.33422342i & 1.03121E-06 & 5.06052E-05 & 15 \\
    13    & 1.16239 & 0.51558364 + 0.33422325i & 1.37117E-06 & 6.92931E-05 & 14 \\
    14    & 1.16731 & 0.89785308 + 0.29205152i & 1.70880E-07 & 9.00285E-05 & 18 \\
\bottomrule
\end{array}
$
\caption{Results obtained using the iterative method \eqref{eq:S3-005}.}\label{tab:01}
\end{table}

\end{example}

\begin{example}

Let $f$ be a function, with

\begin{align*}
f(x)=&-96.98x^{15}-96.82x^{14}-3.87x^{13}+25.78x^{12}\\
&+90.68x^{11}+48.05x^{10}+50.54x^9-5.16x^8\\
&+47.01x^7+90.23x^6+87.09x^5+53.09x^4\\
&+15.38x^3+97.98x^2-61.98x+14.69.
\end{align*}

Then the initial condition $x_0=0.15$ is chosen to use the iterative method  given by \eqref{eq:S3-005}. Consequently, we obtain the results of the Table \ref{tab:02}.

\begin{table}[!ht]
\centering
$
\begin{array}{c|ccccc}
\toprule
&\alpha& x_n&\norm{x_n - x_{n-1} }_2  &\norm{f\left(x_n \right)}_2& n \\ 
\midrule
    1     & 0.88451 & 0.7739975 - 0.54762173i & 1.17047E-07 & 5.48998E-05 & 28 \\
    2     & 0.90499 & -0.82526288 + 0.64969528i & 1.14018E-07 & 7.66665E-05 & 24 \\
    3     & 0.90731 & 0.03271742 + 1.02608471i & 1.14018E-07 & 4.71095E-05 & 19 \\
    4     & 0.90863 & -0.48361539 + 0.928383i & 1.16619E-07 & 6.49343E-05 & 28 \\
    5     & 0.90923 & 0.03271738 - 1.02608473i & 1.30384E-07 & 5.48149E-05 & 22 \\
    6     & 0.93627 & 0.28667103 + 0.19437684i & 4.73466E-06 & 3.81068E-05 & 12 \\
    7     & 0.94059 & 0.30392964 + 0.77330882i & 9.96092E-07 & 6.75516E-05 & 16 \\
    8     & 0.94155 & 0.77399751 + 0.54762167i & 1.14018E-07 & 3.18386E-05 & 14 \\
    9     & 0.95179 & 0.2866711 - 0.19437464i & 1.46720E-05 & 9.07085E-05 & 9 \\
    10    & 0.95499 & -1.16959779 + 0.06354745i & 1.01980E-07 & 1.49833E-05 & 15 \\
    11    & 0.99283 & 1.16397068 & 1.90000E-07 & 6.26264E-05 & 6 \\
    12    & 1.04799 & -0.48361536 - 0.92838301i & 2.15407E-07 & 7.67702E-05 & 31 \\
    13    & 1.05463 & 0.30392985 - 0.77330891i & 7.60000E-07 & 5.64786E-05 & 18 \\
    14    & 1.06431 & -1.16959776 - 0.06354748i & 2.00998E-07 & 6.59751E-05 & 14 \\
    15    & 1.06455 & -0.82526289 - 0.64969531i & 4.12311E-08 & 2.72233E-05 & 21 \\
\bottomrule
\end{array}
$
\caption{Results obtained using the iterative method \eqref{eq:S3-005}.}\label{tab:02}
\end{table}

\end{example}

\begin{example}

Let $f$ be a function, with

\begin{align*}
f(x)=&-57.62x^{16}-56.69x^{15}-37.39x^{14}-19.91x^{13}
+35.83^{12}\\
&-72.47^{11}+44.41x^{10}+43.53x^9+59.93x^8\\
&-42.9x^7-54.24x^6+72.12x^5-22.92x^4\\
&+56.39x^3+15.8x^2+60.05x+55.31.
\end{align*}

Then the initial condition $x_0=0.83$ is chosen to use the iterative method  given by \eqref{eq:S3-005}. Consequently, we obtain the results of the Table \ref{tab:03}.

\begin{table}[!ht]
\centering
$
\begin{array}{c|ccccc}
\toprule
&\alpha& x_n&\norm{x_n - x_{n-1} }_2  &\norm{f\left(x_n \right)}_2& n \\ 
\midrule
    1     & 0.81691 & 0.8812118 - 0.42696217i & 1.08167E-07 & 8.16395E-05 & 49 \\
    2     & 0.83851 & 1.03423973 & 6.00000E-08 & 7.69988E-05 & 56 \\
    3     & 0.97383 & -1.0013396 & 6.37436E-05 & 6.64131E-05 & 7 \\
    4     & 0.99055 & -0.35983764 + 1.18135267i & 6.70820E-08 & 2.53547E-05 & 21 \\
    5     & 0.99059 & -0.70050491 + 0.78577099i & 1.70294E-07 & 9.13799E-06 & 17 \\
    6     & 0.99219 & -0.70050494 - 0.78577099i & 1.18229E-06 & 5.28258E-05 & 7 \\
    7     & 0.99283 & 0.36452491 - 0.83287828i & 3.63610E-06 & 4.30167E-05 & 17 \\
    8     & 0.99347 & -0.28661378 - 0.8084062i & 1.38226E-05 & 9.04752E-05 & 8 \\
    9     & 0.99427 & -0.35983765 - 1.18135267i & 1.26491E-07 & 4.09162E-05 & 16 \\
    10    & 0.99539 & -1.3699527 & 2.30000E-07 & 7.02720E-05 & 14 \\
    11    & 1.12775 & -0.62435238 & 1.25000E-06 & 6.46233E-05 & 4 \\
    12    & 1.16423 & 0.58999229 - 0.86699687i & 7.07107E-08 & 7.38972E-05 & 25 \\
    13    & 1.16595 & 0.36452487 + 0.83287805i & 3.06105E-07 & 9.42729E-05 & 15 \\
    14    & 1.16607 & 0.58999222 + 0.86699689i & 5.00000E-08 & 5.09054E-05 & 18 \\
    15    & 1.16647 & 0.88121183 + 0.42696223i & 4.12311E-08 & 5.37070E-05 & 39 \\
    16    & 1.20923 & -0.28661363 + 0.8084063i & 1.94165E-07 & 5.02799E-05 & 16 \\
\bottomrule
\end{array}
$
\caption{Results obtained using the iterative method \eqref{eq:S3-005}.}\label{tab:03}
\end{table}

\end{example}

\section{Conclusions}

The F N-R method is very efficient to find roots of polynomials since it does not present the divergence problems, like the classical N-R method, for a polynomial with only complex roots when using real initial conditions. However, the really interesting thing is that this method opens up the possibility of creating new fractional iterative methods  \cite{torreshern2020,gdawiec2019visual,gdawiec2020newton,akgul2019fractional,cordero2019variant,torres2020fractional}, as well as opens the possibility of creating new hybrid iterative methods by combining the F N-R method with existing iterative methods \cite{plato2003concise,stoer2013}. So in this work it has been given one more application to fractional calculus and has opened the possibility of extending the capacity of the iterative methods that allow us to find zeros of functions more general than polynomials \cite{torres2020reduction,torreshern2020}.

\bibliography{Biblio}

\begin{thebibliography}{10}

\bibitem{plato2003concise}
Robert Plato.
\newblock {\em Concise numerical mathematics}.
\newblock Number~57. American Mathematical Soc., 2003.

\bibitem{stoer2013}
Josef Stoer and Roland Bulirsch.
\newblock {\em Introduction to numerical analysis}, volume~12.
\newblock Springer Science \& Business Media, 2013.

\bibitem{brambila2017fractal}
Fernando Brambila.
\newblock {\em Fractal Analysis: Applications in Physics, Engineering and
  Technology}.
\newblock IntechOpen, 2017.

\bibitem{martinez2017applications2}
Benito~F. Mart{\'\i}nez-Salgado, Rolando Rosas-Sampayo, Anthony
  Torres-Hern{\'a}ndez, and Carlos Fuentes.
\newblock Application of fractional calculus to oil industry.
\newblock {\em Fractal Analysis: Applications in Physics, Engineering and
  Technology}, 2017.
\newblock
  \url{https://www.intechopen.com/books/fractal-analysis-applications-in-physics-engineering-and-technology}.

\bibitem{miller93}
Kenneth~S. Miller and Bertram Ross.
\newblock {\em An introduction to the fractional calculus and fractional
  differential equations}.
\newblock Wiley-Interscience, 1993.

\bibitem{oldham74}
Keith Oldham and Jerome Spanier.
\newblock {\em The fractional calculus theory and applications of
  differentiation and integration to arbitrary order}, volume 111.
\newblock Elsevier, 1974.

\bibitem{hilfer00}
Rudolf Hilfer.
\newblock {\em Applications of fractional calculus in physics}.
\newblock World Scientific, 2000.

\bibitem{torres2020reduction}
A.~Torres-Hernandez, F.~Brambila-Paz, P.M. Rodrigo, and E.~{De-la}-Vega.
\newblock Reduction of a nonlinear system and its numerical solution using a
  fractional iterative method.
\newblock {\em Journal of Mathematics and Statistical Science}, 2020.
\newblock
  \url{http://www.ss-pub.org/wp-content/uploads/2020/10/JMSS2020070201.pdf}.

\bibitem{torres2020nonlinear}
A.~Torres-Hernandez, F.~Brambila-Paz, and J.J. Brambila.
\newblock A nonlinear system related to investment under uncertainty solved
  using the fractional pseudo-newton method.
\newblock {\em Journal of Mathematical Sciences: Advances and Application},
  2020.
\newblock
  \url{http://scientificadvances.co.in/admin/img_data/1470/images/JMSAA7100122150ATorresHernandez.pdf}.

\bibitem{kilbas2006theory}
A.A. Kilbas, H.M. Srivastava, and J.J. Trujillo.
\newblock {\em Theory and Applications of Fractional Differential Equations}.
\newblock Elsevier, 2006.

\bibitem{torreshern2020}
A.~Torres-Hernandez, F.~Brambila-Paz, and E.~{De-la}-Vega.
\newblock Fractional newton-raphson method and some variants for the solution
  of nonlinear systems.
\newblock {\em Applied Mathematics and Sciences: An International Journal
  (MathSJ)}, 2020.
\newblock \url{https://airccse.com/mathsj/papers/7120mathsj02.pdf}.

\bibitem{gdawiec2019visual}
Krzysztof Gdawiec, Wies{\l}aw Kotarski, and Agnieszka Lisowska.
\newblock Visual analysis of the newton’s method with fractional order
  derivatives.
\newblock {\em Symmetry}, 11(9):1143, 2019.

\bibitem{gdawiec2020newton}
Krzysztof Gdawiec, Wies{\l}aw Kotarski, and Agnieszka Lisowska.
\newblock Newton’s method with fractional derivatives and various iteration
  processes via visual analysis.
\newblock {\em Numerical Algorithms}, pages 1--58, 2020.

\bibitem{akgul2019fractional}
Ali Akg{\"u}l, Alicia Cordero, and Juan~R Torregrosa.
\newblock A fractional newton method with 2$\alpha$th-order of convergence and
  its stability.
\newblock {\em Applied Mathematics Letters}, 98:344--351, 2019.

\bibitem{cordero2019variant}
Alicia Cordero, Ivan Girona, and Juan~R Torregrosa.
\newblock A variant of chebyshev’s method with 3$\alpha$th-order of
  convergence by using fractional derivatives.
\newblock {\em Symmetry}, 11(8):1017, 2019.

\bibitem{torres2020fractional}
A.~Torres-Hernandez, F.~Brambila-Paz, P.M. Rodrigo, and E.~{De-la}-Vega.
\newblock Fractional pseudo-newton method and its use in the solution of a
  nonlinear system that allows the construction of a hybrid solar receiver.
\newblock {\em Applied Mathematics and Sciences: An International Journal
  (MathSJ)}, 2020.
\newblock \url{https://airccse.com/mathsj/papers/7220mathsj01.pdf}.

\end{thebibliography}
\bibliographystyle{unsrt}

\end{document}